\documentclass[12pt]{article}
\usepackage{amsmath,amsthm,amsfonts,latexsym,amsopn,verbatim,amscd,amssymb}
\theoremstyle{plain}
\newtheorem{Thm}{Theorem}[section]
\newtheorem{Lem}[Thm]{Lemma}

\theoremstyle{definition}

\newtheorem{Def}[Thm]{Definition}

\numberwithin{equation}{section}

\DeclareMathOperator{\idem}{idem}

\DeclareMathOperator{\U}{U}

\DeclareMathOperator{\win}{Win}

\DeclareMathOperator{\In}{In}
\DeclareMathOperator{\J}{J}

\newcommand{\bnum}{\begin{enumerate}}
\newcommand{\enum}{\end{enumerate}}
\begin{document}
\begin{center}
\textbf{WEAK CLEAN INDEX OF A RING }\\
\end{center}
\begin{center}
Dhiren Kumar Basnet\\
\small{\it Department of Mathematical Sciences, Tezpur University,
 \\ Napaam, Tezpur-784028, Assam, India.\\
Email: dbasnet@tezu.ernet.in}
\end{center}

\begin{center}
Jayanta Bhattacharyya \\
\small{\it Department of Mathematical Sciences, Tezpur University,
 \\ Napaam, Tezpur-784028, Assam, India.\\
Email: jbhatta@tezu.ernet.in}
\end{center}
\noindent \textit{\small{\textbf{Abstract:}  }} Motivated by the concept of clean index of  rings \cite{cir}, we introduce the concept of weak clean index of rings. For any element $a$ of a ring $R$ with unity, we define $ \chi(a)=\{e\in R\mid e^2=e\text{ and }a-e \mbox{ or } a+e \mbox{ is a unit}\}$. The weak clean index of $R$ is defined as $\sup \{|\chi(a)|: a\in R\}$ and it is denoted by $\win(R),$ where $| \chi(a)| $ denotes the cardinality of the set $\chi(a)$. In this article, we characterize rings of weak clean indices $1$, $2$ and $3$.
\bigskip

\noindent \small{\textbf{\textit{Key words:}} Clean ring, weak clean ring, clean index, weak clean index.} \\
\smallskip

\noindent \small{\textbf{\textit{$2010$ Mathematics Subject Classification:}}  16N40, 16U99.} \\
\smallskip

\bigskip

\section{INTRODUCTION}
In this article we assume ring $R$ to be associative with unity unless otherwise stated, we also assume modules (and bimodules) to be unitary. The Jacobson radical, group of units, set of idempotents and set of nilpotent elements of a ring $R$ are denoted by J$(R)$, U$(R)$, idem$(R)$ and nil$(R)$ respectively. Notion of \textit{clean ring} was first introduced by Nicholson \cite{wkn}, which was later generalized by  Ahn and Andreson\cite{and} as \textit{weakly clean ring} defined as a ring $R$ in which  each element $r\in R$ can be written as $r=u+e$ or $r=u-e$ for some $u\in \U(R)$ and $e\in \idem(R)$.  Further, Lee and Zhou \cite{cir, cir2} introduced and studied \textit{clean index }of rings, which actually motivated us to introduce and study weak clean index of rings.
\begin{Def}
For any element $a$ of $R,$ we define $ \chi(a)=\{e \in \idem(R) \mid a-e\in \U(R) \mbox{ or } a+e\in\U(R)\}$. The weak clean index of $R$ is defined as $\sup \{|\chi(a)|: a\in R\}$ and it is denoted by $\win(R),$ where $| \chi(a)| $ denotes the cardinality of the set $\chi(a)$.
\end{Def}

\section{BASIC PROPERTIES}

Some basic properties related to weak clean index are presented here as a preparation for the article.
\begin{Lem}\label{5} Let $R$ be a ring and $e, a, b\in R$. Then following hold:
\begin{enumerate}
  \item  For a central nilpotent $n\in R$, $|\chi(n)|=1$. Whereas for a central idempotent $e \in R$, $|\chi(e)| \geq1$, thus for any ring $\win(R)\geq1$.
  \item If $a-b\in\J(R)$ then $| \chi(a)|=|\chi(b)|$.
  \item $e\in\chi(a)$ then $1-e\in \chi(1-a) $ or $1-e\in \chi(1+a) $. Converse holds if $2\in\J(R)$.
  \item Let $\sigma$ be an automorphism or anti-automorphism of $R$. Then $e\in\chi(a)$ iff $\sigma(e)\in \chi(\sigma(a))$; so $| \chi(a)|=|\chi(\sigma(a))|$. In particular $|\chi(a)|=|\chi(uau^{-1})|$, where $u$ is a unit of $R$.
  \item If a ring $R$ has at most $n$ units or at most $n$ idempotents, then $\win(R)\leq n$. In particular, if $R$ is a local ring then $\win(R)\leq 2$.
  \item If $R$ is local, then $\win(R)=2$ iff $R/\J(R)\ncong \mathbb{Z}_2$.
  \item Let $R$ be a clean ring with $2\in  \U(R)$. Then $\win(R) = | \chi(2^{-1}) |$, in other words $\idem(R)=\chi(2^{-1})$.
\end{enumerate}
\end{Lem}
$Proof$ (i) Let $a$ be a central nilpotent such that $a^n=0$ for some $n\in \mathbb{N}$, then we have  $a=(a+1)-1$ is a weak clean expression, hence $1\in \chi(a)$; $| \chi(a)|\geq 1$. If possible let $e\neq 1\in \chi(a)$, then there exists a $u\in\U(R)$ such that $a=u+e\mbox{ or }u-e$. If $a=u-e$, by using binomial expansion and the fact that $a^n=0$ we have
$0=(u-e)^n=u^n-{n \choose 1}eu^{n-1}+{n \choose 2}eu^{n-2}-\dots +(-1)^{n-1}eu+(-1)^ne$, implies $u^n\in eR$ contradicting the fact that $e\neq 1$. Next if $a=u+e$, similarly we get a contradiction.  let $e$ be a central idempotent then $e=1-(1-e)$ i.e., $ 1-e \in \chi(e)$ and  therefore $|\chi(e)|\geq1$. \\
(ii) Let $w=a-b\in\J(R)$. If $e\in\chi(a)$, we have $a+e\in \U(R)$ or $a-e\in \U(R)$.\\
Case 1: $u=a+e\in\ U(R) $ implies $u=b+w+e \Rightarrow b+e=u-w \in \U(R)$ and therefore $e\in\chi(b)$.\\
Case 2: $v=a-e\in \U(R)$  implies $b-e=v-w\in \U(R)$ and therefore $e\in\chi(b)$.\\
Thus $\chi(a)\subseteq \chi(b)$ and  by symmetry we get  $\chi(b)\subseteq \chi(a)$ resulting $\chi(a)=\chi(b)$.\\
(iii)  Let $e\in \chi(a)$, then we have $a+e\in \U(R)$ or $a-e\in \U(R)$. Let $a-e\in\U(R)$ then we have $(1-a)-(1-e)=e-a \in \U(R)$, so $1-e\in\chi(1-a)$. Similarly if $a+e\in \U(R)$ then we have $(1+a)-(1-e)=a+e \in \U(R)$, therefore $1-e\in\chi(1+a)$. Conversely if $(1-e)\in\chi(1-a)$, we have $(1-a)- (1-e)=u\in\U(R)$ or $(a-1)+(1-e)=v\in\U(R)$, that is $a-e=-u$ or $a-e=v$, so in this case $e\in\chi(a)$. If $(1-e)\in\chi(1+a)$, we have $(1+a)- (1-e)=u\in\U(R)$ or $(a+1)+(1-e)=v\in\U(R)$ this implies $a+e=u$ or $a-e=v-2 \in\U(R)$, as $2\in\J(R)$, hence we get $e\in\chi(a)$.\\
(iv)and (v) are straightforward.\\
(vi) If $R$ is a local ring then we have $\win(R)\leq 2$, as $\idem(R)=\{0,1\}$. Let $R/\J(R)\cong \mathbb{Z}_2$,  that is, $R$ is uniquely clean. If possible let $\win(R)=2$, that is there exists at least an element $a\in R$ such that $\{0,1\}=\chi(a)$. So, $a\in\U(R)$ and $a-1 \in \U(R)$ or $a+1 \in \U(R)$. If $a\in\U(R)$ and $u=a-1 \in \U(R)$, then we have two clean expressions for $a$ which is a contradiction. Similarly $a\in\U(R)$ and $u=a+1 \in \U(R)$, give two  clean expressions for $u$, which is a contradiction, hence $\win(R)\neq 2$. Conversely, let $\win(R)=1$, then $\In(R)=1$ as $\In(R)\leq \win(R)$, hence the result follows by \textbf{Theorem $2.1$} of \cite{hc1}\\
(vii) Let $e\in \idem(R)$ and let $2 \in \U(R)$, now we have $(2^{-1}-e)\in\U(R)$, as $2(1-2e)$ is inverse of $2^{-1}-e$, therefore $\idem(R)\subseteq \chi(2^{-1})$, hence $\win(R) = | \chi(2^{-1}) |$. $\square$\\

In a ring $R$, $q\in R$ is called quasi-regular element, if there is a $p\in R$, such that
$q+p+qp = 0 = p+q+pq$. The set of all all quasi-regular elements of ring $R$ is denoted by $Q(R)$.
\begin{Lem}\label{winlm0}
  If $S$ is a subring of a ring $R$, where $R$ and $S$ may not share same identity, then $\win(S)\leq\win(R)$.
\end{Lem}
$Proof$ For $a\in R$, let $J(a)=J_1(a)\cup J_2(a)$, where $J_1(a)=\{ q\in Q(R) : (a-q)^2=a-q\}$ and $J_2(a)=\{ q\in Q(R) : (q-a)^2=q-a\}$. Claim: $\win(R)=\sup\{|J(b)|: b\in R\}$. Note that $\U(R)=\{ 1+q : q\in Q(R)\}$. For any $a\in R$, $\chi(a)=\{ (a-1)-j : j\in J_1(a-1)\}\cup \{j-(a-1) : j\in \J_2(a-1)\}$, therefore $|\chi(a)|=|J(a-1)|$. Thus $\win(R)=\sup\{J(b) : b\in R\}$. Because $ Q(S)\subseteq Q(R) $ it follows that $\win(S)\leq \win(R)$.\\
\textbf{Proof of $|\chi(a)|=|J(a-1)|$:}\\
Let $e\in\chi(a)$, we have
\begin{align*}
  &\Leftrightarrow a-e=u \mbox{ or }  a+e=u, \mbox{ for some } u\in \U(R). \\
  &\Leftrightarrow a-u=e \mbox{ or } u-a=e,\\
  &\Leftrightarrow a-1-q=e \mbox{ or } 1+q-a=e ,\mbox{ for some } q=1+u \mbox{ as } \U(R)=1+Q(R).\\
  &\Leftrightarrow (a-1)-q=e \mbox{ or }  q-(a-1)=e, \\
  &\Leftrightarrow e\in J_1(a-1) \mbox{ or } e\in J_2(a-1). \square
\end{align*} 
\begin{Thm}
  Let $k\geq 1$ be an integer, then the following are equivalent for a ring $R$
\begin{enumerate}
  \item $\win(R[[x]])=k$.
  \item $\win(R[x])=k$.
  \item  $R$ is abelian and $\win(R)=k$.
\end{enumerate}
\end{Thm}
$Proof$ By above \textbf{Lemma \ref{winlm0}} we have $\win(R)\leq \win(R[x]) \leq \win(R[[x]])$. Suppose that $R$ is not abelian and $e$ be a non-central idempotent of $R$, let $er\neq re$ for some $r\in R$. So either $er(1-e)\neq 0$ or $(1-e)re\neq 0$, without loss of generality we may assume that $er(1-e)\neq 0$. For $i=1, 2, 3, \dots $
\begin{align*}
  a:= & [1+er(1-e)]-e \\
  = & [1+er(1-e)(1+x^i)]-[e+er(1-e)x^i]
\end{align*}
are infinitely many distinct weak clean expressions of $a$ in $R[x]$. Now suppose $R$ is abelian, it is easy to see that idempotents of $R[[x]]$ are all in $R$ and for any $\alpha=a_0+a_1x+a_2x^2+\dotsm \in R[[x]]$,  $\chi_{R[[x]]}(\alpha)\subseteq \chi_R(a_0).$ Thus $|\chi(\alpha)| \leq |\chi(a_0)|$, therefore $\win(R[[x]])\leq \win(R)$ and consequently the result follows. $\square$
\section{Rings with weak clean index 1, 2 and 3}
In this section we try to characterize the rings of weak clean index 1, 2 and 3.
\begin{Thm}\label{winth1}
$\win(R)=1$ iff $R$ is abelian and for any $0\neq e^2=e\in R$, $e\neq u+v$ for any $u, v\in\U(R)$.
\end{Thm}
$Proof$ $(\Rightarrow)$ Let $e^2=e\in R$. For any $r\in R$, $1-e=[1+ er(1-e)] -[e+er(1-e)]$ are two weak clean expression of $1-e$; so $e=[e+er(1-e)]$, that is $re=ere$, similarly we have $er=ere$, so $R$ is abelian. Suppose that $0\neq e^2=e\in R$, $e = u+v$ for some $u, v\in\U(R)$, then $v=v+0=-u+e$ are two weak clean expressions of $v$, implies $| \chi(v) |= 2 $, which is a contradiction.\\
$(\Leftarrow)$ Let $a\in R$ has two weak clean expressions, $a = u_1 + e_1$ or $u_1 - e_1$ and
$a = u_2 + e_2$ or $u_2 - e_2$, for $e_1,e_2 \in \idem(R)$, $e_1\neq e_2$ and $u_1,u_2\in\U(R)$.\\
\textbf{Case I:} if $a=u_1+e_1=u_2+e_2$, we have $e_1-e_2=u_2-u_1.$ Define $f:=e_1(1-e_2)$, then $f=f^2\in R$. Now
\begin{align*}
  f= & [e_2+(u_2-u_1)](1-e_2) \\
  = & u_2(1-e_2)-u_1(1-e_2)\\
  = & [u_2(1-e_2)+e_2]-[u_1(1-e_2)+e_2].
\end{align*}
As $[u_2(1-e_2)+e_2], [u_1(1-e_2)+e_2]$ are units in $R$, 
therefore $f=0$, hence $e_1=e_1e_2$. Similarly we have $e_2 = e_1e_2$,  a contradiction, thus $\chi(a)\leq 1$.\\
\textbf{Case II:} If $a=u_1+e_1=u_2-e_2$, implies $e_1+e_2=u_2-u_1.$ Define $f:=e_1(1-e_2)$, then $f=f^2 \in R$. We have
\begin{align*}
  f= & [-e_2+(u_2-u_1)](1-e_2) \\
  = & u_2(1-e_2)-u_1(1-e_2)\\
  = & [u_2(1-e_2)+e_2]-[u_1(1-e_2)+e_2].
\end{align*}
As $[u_2(1-e_2)+e_2], [u_1(1-e_2)+e_2]$ are units in $R$, 
so $f=0$, hence $e_1=e_1e_2$. Similarly we get $e_2=e_1e_2$, that is $e-1=e_2$,  a contradiction, thus $\chi(a)\leq 1$.\\
\textbf{Case III:} If $a=u_1-e_1=u_2-e_2$, we have $e_1-e_2=u_1-u_2.$ Define $f:=e_1(1-e_2)$, then $f=f^2 \in R$.
\begin{align*}
  f= & [e_2+(u_1-u_2)](1-e_2) \\
  = & u_1(1-e_2)-u_2(1-e_2)\\
  = & [u_1(1-e_2)+e_2]-[u_2(1-e_2)+e_2].
\end{align*}
Since $R$ is abelian, so $[u_2(1-e_2)+e_2], [u_1(1-e_2)+e_2]$ are units. 
So, as in above case $\chi(a)\leq 1$. Thus combining above cases we conclude that $\win(R)=1$.$\square$
\begin{Lem}\label{winla1}
  Let $R=A\times B$ be direct product of rings $A$ and $B$, such that $\win(A)=1$, then $\win(R)=\win(B)$.
\end{Lem}
$Proof$ Since $A, B$ are subrings of $R$, so by \textbf{Lemma \ref{winlm0}}, $\win(B)\leq \win(R)$. If $\win(B)=\infty$, then $\win(R)=\infty$ and so we have $\win(R)=\win(B)$. Next let $\win(B)=k<\infty$ be a positive integer, so there is a $b\in B$, such that $|\chi(b)|=k$. Now for $(0,b)\in R$, $|\chi(0,b)|=k$, hence $\win(R)\geq k$. Suppose that $\win(R)>k$. Then there exists $(a,b)\in R$ such that $(a,b)$ has at least $k+1$ weak clean expressions in $R$. Let $g$ be an integer such that $1\leq g\leq k$  and let
$(a,b)=\left\{
         \begin{array}{ll}
           (u_i,v_i) + (e_i,f_i), & \hbox{i= 1, 2, 3, \dots, g} \\
           (u_j,v_j) - (e_j,f_j), & \hbox{j= g+1, g+2, \dots, k, k+1.}
         \end{array}
       \right.$ \\
        are $k+1$ distinct weak clean expressions for $(a,b)$, such that no two $(e_i,f_i)$'s are equal. Now, $a=u_i + e_i=u_j - e_j$, $(i= 1, 2, 3, \dots, g \mbox{ and } j= g+1, g+2, \dots, k+1)$ are weak clean expressions of $a$ in $S$. Since $|\chi(a)| \leq 1$, so all $e_i's$ and $e_j's$ are equal. So 
 $k+1=|\chi((a,b))| \\
 = |\{(e_i,f_i), (e_j,f_j)| i= 1, 2, \dots, g \mbox{ and } j= g+1, g+2, \dots, k+1\}| \\
  =  |\{e_i,e_j| i= 1, 2, \dots, g; j= g+1,  \dots, k\}|\times |\{f_i, f_j| i= 1, 2, \dots, g;j= g+1,  \dots, k\}| \\
  =  |\chi(a)|\times |\chi(b)|\\
  =  |\chi(b)|$

\noindent which is a contradiction and this completes the proof. $\square$

\begin{Def}
  A ring $R$ is said to be elemental if idempotents of $R$ are trivial and $1=u+v$ for some $u, v \in \U(R)$.
\end{Def}

\begin{Thm}
  For a ring $R$, $\win(R)=2$ iff one of the following holds:
\begin{enumerate}\label{winth2}
  \item $R$ is elemental.
  \item $R = A \times B$, where $A$ is elemental ring and $\win(B)=1$.
  \item $R=\left(
                 \begin{array}{cc}
                   A & M \\
                   0 & B \\
                 \end{array}
               \right)$, where $\win(A)=\win(B)=1$ and $_AM_B$is a bimodule with $|M|=2$.
\end{enumerate}
\end{Thm}
$Proof$ $(\Leftarrow)$ If (i) holds then by the definition of elemental ring, we have $1=u+v$ for some $u, v \in\U(R)$, therefore by \textbf{Theorem \ref{winth1}}, $\win(R)>1$. By \textbf{Lemma \ref{winlm0}} $\win(R)\leq |\idem(R)|=2$, thus $\win(R)=2$.\\
If (ii) holds then $\win(R)=2$ by \textbf{(i)} and \textbf{Lemma \ref{winla1}}.\\
If (iii) holds, for $\alpha_0=\left(\begin{array}{cc}
                                0 & 0 \\
                                0 & 1
                              \end{array}\right),$ we have
$\left\{ \left(\begin{array}{cc}
          1 & w \\
          0 & 0
        \end{array}\right) : w\in M\right\}\subseteq \chi(\alpha_0)$. So, $\win(R)\geq |\chi(\alpha_0)|\geq |M|=2$. For any $\alpha=\left(\begin{array}{cc}
                                a & x \\
                                0 & b
                              \end{array}\right) \in R$, $$|\chi(\alpha)|=\left| \left\{ \left(\begin{array}{cc}
                                e & w \\
                                0 & f
                              \end{array}\right)\in R : e\in \chi(a), f\in\chi(b), w=ew+wf \right\}\right|.$$
Since $|M|=2,$ $ |\chi(a)| \leq 1$ and $|\chi(b)|\leq 1$, it follows that $|\chi(\alpha)|\leq 2$ and hence $\win(R)=2$.\\
$(\Rightarrow)$ Suppose $R$ is abelian, as $\win(R)\neq 1$, so there exists $(0\neq) e=e^2\in R$ such that $e=u+v$, where $u, v\in \U(R)$. So we have $e=eu+ev$, where $eu, ev\in \U(eR)$, hence $\win(eR)\leq 2$. But $\win(eR)\leq \win(R) = 2$ by \textbf{Lemma \ref{winlm0}}. From $R = A \times B$ where $A=eR$ and $B=(1-e)R$ it follows, $\win(B)=1$. If $A$ has a non trivial idempotent $f$ then $A = fA + (1-f)A$ where, $f = fu + fv$ and $ 1-f = (1-f)u+(1-f)v $, where $ fu, fv \in \U(fA)$ and $ (1-f)u, (1-f)v \in \U((1-f)A) $, thus by \textbf{Theorem 5} \cite{cir}, we have $\In(fA)\geq 2$ and $\In((1-f)A)\geq 2$, so $\In(A) \geq 2\times 2=4$, as $\In(R) \leq \win(R)$, hence a contradiction. Thus (i) holds if $e=1$ and (ii) holds if $e\neq 1$. Suppose $R$ is not abelian and let $e^2= e \in R$ be a non central idempotent. If neither of $eR(1-e)$ and $(1-e)Re$ is zero, then take $0\neq x\in eR(1-e)$ and $0\neq y\in (1-e)Re$ to get
$ 1-e = (1+x)-(x+e) = (1+y)-(y+e)$. Therefore $|\chi(1-e)|\geq 3$, which is a contradiction. So without loss of generality we can assume that  $eR(1-e) \neq 0$ and $(1-e)Re=0$. The Peirce decomposition of $R$ gives $$R=\left(
                                     \begin{array}{cc}
                                       eRe & eR(1-e) \\
                                       0 & (1-e)R(1-e) \\
                                     \end{array}
                                   \right).$$
As above $2=\win(R)\geq |eR(1-e)|$; so $|eR(1-e)| = 2$. Write $eR(1-e)=\{0, x\}$. If $\win(eRe)=2$, then there exists an $a\in R$ such that $|\chi(a)| = 2$. Thus we have following cases\\
\textbf{Case I:} Let $ a= u_1 + e_1 = u_2 + e_2 $, where $ u_1, u_2 \in \U(eRe) $ and $ e_1, e_2 \in \idem(eRe) $. If $ e_1x =0 $, we have for $ A=\left( \begin{array}{cc}   a & 0 \\    0 & 0 \\  \end{array} \right) \in R $\\
$
A = \left( \begin{array}{cc}   u_1 & 0 \\    0 & -1 \\  \end{array} \right) +
\left( \begin{array}{cc}   e_1 & 0 \\    0 & 1 \\  \end{array} \right) =
\left( \begin{array}{cc}   u_2 & 0 \\    0 & -1 \\  \end{array} \right) +
\left( \begin{array}{cc}   e_2 & 0 \\    0 & 1 \\  \end{array} \right) =
\left( \begin{array}{cc}   u_1 & x \\    0 & -1 \\  \end{array} \right) +
\left( \begin{array}{cc}   e_1 & x \\    0 & 1 \\  \end{array} \right)
$
are three distinct weak clean expressions of $A$ in $R$, which implies $|\chi(A)| \geq 3$, that is a contradiction. If $ex=x$ then for $B=\left( \begin{array}{cc}   a & 0 \\    0 & 1 \\  \end{array} \right)$, we have\\
$
B = \left( \begin{array}{cc}   u_1 & 0 \\    0 & 1 \\  \end{array} \right) +
\left( \begin{array}{cc}   e_1 & 0 \\    0 & 0 \\  \end{array} \right) =
\left( \begin{array}{cc}   u_2 & 0 \\    0 & 1 \\  \end{array} \right) +
\left( \begin{array}{cc}   e_2 & 0 \\    0 & 0 \\  \end{array} \right) =
\left( \begin{array}{cc}   u_1 & x \\    0 & 1 \\  \end{array} \right) +
\left( \begin{array}{cc}   e_1 & x \\    0 & 0 \\  \end{array} \right)
$
are three distinct weak clean expressions of $B$ in $R$, which implies $|\chi(B)| \geq 3$, that is a contradiction.\\
\textbf{Case II:} Let $ a = u_1 - e_1 = u_2 + e_2 $,  where $ u_1, u_2 \in \U(eRe) $ and $ e_1, e_2 \in \idem(eRe) $. So if $ e_1x =0 $, we have for $ A=\left( \begin{array}{cc}   a & 0 \\    0 & 0 \\  \end{array} \right) \in R $\\
$
A = \left( \begin{array}{cc}   u_1 & 0 \\    0 & 1 \\  \end{array} \right) -
\left( \begin{array}{cc}   e_1 & 0 \\    0 & 1 \\  \end{array} \right) =
\left( \begin{array}{cc}   u_2 & 0 \\    0 & -1 \\  \end{array} \right) +
\left( \begin{array}{cc}   e_2 & 0 \\    0 & 1 \\  \end{array} \right) =
\left( \begin{array}{cc}   u_1 & x \\    0 & 1 \\  \end{array} \right) -
\left( \begin{array}{cc}   e_1 & x \\    0 & 1 \\  \end{array} \right)
$
are three distinct weak clean expressions of $B$ in $R$, which implies $|\chi(B)| \geq 3$, that is a contradiction. If $ex=x$ then for $B=\left( \begin{array}{cc}   a & 0 \\    0 & 1 \\  \end{array} \right)$, we have\\
$
B = \left( \begin{array}{cc}   u_1 & 0 \\    0 & 1 \\  \end{array} \right) -
\left( \begin{array}{cc}   e_1 & 0 \\    0 & 0 \\  \end{array} \right) =
\left( \begin{array}{cc}   u_2 & 0 \\    0 & 1 \\  \end{array} \right) +
\left( \begin{array}{cc}   e_2 & 0 \\    0 & 0 \\  \end{array} \right) =
\left( \begin{array}{cc}   u_1 & x \\    0 & 1 \\  \end{array} \right) -
\left( \begin{array}{cc}   e_1 & x \\    0 & 0 \\  \end{array} \right)
$
are three distinct weak clean expressions of $B$ in $R$, which implies $|\chi(B) | \geq 3$, again a contradiction.\\
\textbf{Case III:} Let $ a = u_1 - e_1 = u_2 - e_2 $,  where $ u_1, u_2 \in \U(eRe) $ and $ e_1, e_2 \in \idem(eRe) $, we get a contradiction similar to case I.\\
This shows that $\win(eRe)=1$, similarly $\win((1-e)R(1-e))$ =1.
\begin{Thm}
  $\win(R)=3$ iff $ R = \left( \begin{array}{cc}   A & M \\    0 & B \\  \end{array} \right)$, where $\win(A)=\win(B)=1$ and $_AM_B$is a bimodule with $|M|=3$.
\end{Thm}
$Proof$ $(\Leftarrow)$ For $\alpha_0=\left(\begin{array}{cc} 0 & 0 \\ 0 & 1 \end{array}\right)$, we have
$\left\{ \left(\begin{array}{cc} 1 & w \\ 0 & 0 \end{array}\right) : w\in M\right\}\subseteq \chi(\alpha_0)$. So, $\win(R)\geq |\chi(\alpha_0)|\geq |M|=3$. For any $\alpha=\left(\begin{array}{cc} a & x \\ 0 & b\end{array}\right) \in R$,

 $$|\chi(\alpha)|=\left| \left\{ \left(\begin{array}{cc}e & w \\ 0 & f \end{array}\right)
       \in R : e\in \chi(a), f\in\chi(b), w=ew+wf \right\}\right|,$$
as $|M|=3,$ $|\chi(a)| \leq 1$ and $|\chi(b)|\leq 1$ it follows $|\chi(\alpha)|\leq 3$, hence $\win(R)=3$.

$(\Rightarrow)$ Suppose $\win(R)=3$, from the proof of \textbf{Theorem \ref{winth2}} we see that an abelian ring not satisfying condition (i) and (ii), contains a subring whose weak clean index is greater than $4$, therefore $R$ must be non abelian.\\
Let $e$ be a non central idempotent in the ring $R$, then Peirce decomposition of $R$ gives $$R=\left(\begin{array}{cc} eRe & eR(1-e) \\ (1-e)Re & (1-e)R(1-e)\end{array}\right).$$ Let $A = eRe,$ $B = (1-e)R(1-e),$ $M = eR(1-e)$, $N = (1-e)Re$. Suppose $|M| \neq 0$ and $|N| \neq 0$. As $\chi(1-e) \supseteq \{ e-x, e-y : x\in M, 0 \neq y \in N \}$, it follows that $ 3= \win(R) \geq |\chi(1-e)| > |M| + |N| - 1 $, therefore $ |M| = |N| = 2$. Write
$M = \{ 0, x \}$, $N = \{ 0, y \}$. Note that $2x = 0 = 2y$. If $xyx = 0$, then
$(x + y + xy + yx)^4 = 0$ and $\chi(1-e) \supseteqq \{ e, e-x, e-y, e + x + y + xy + yx \}$,
so $\win(R) \geq 4$, a contradiction. If $yxy = 0$, then $(x + y + xy + yx)^4 = 0$ and $\chi(2-e) \supseteqq \{ 1-e, 1-e+x, 1-e+y, 1-e + x + y + xy + yx \}$,
therefore $\win(R) \geq 4$, a contradiction. Hence $xyx \neq 0$ and $yxy \neq 0$. It follows that $xyx = x$ and $yxy = 0$. Let $f = xy$ and $g = yx$, then clearly $f, g$ are idempotents. So we have $R \supseteq L:= \left(\begin{array}{cc} fRf & M \\ N & gRg\end{array}\right).$ By \textbf{Lemma \ref{winlm0}}, $\win(L) \leq 3$, but for
$\alpha = \left(\begin{array}{cc}0 & x \\ y & g \end{array}\right)$ we have
\begin{align*}
  \alpha= & \left(\begin{array}{cc}0 & x \\ y & 0 \end{array}\right)  +
              \left(\begin{array}{cc}0 & 0 \\ 0 & g \end{array}\right)\\
  = & \left(\begin{array}{cc}0 & x \\ y & g \end{array}\right)  +
              \left(\begin{array}{cc}0 & 0 \\ 0 & 0 \end{array}\right)\\
  = & \left(\begin{array}{cc}f & x \\ y & 0 \end{array}\right)  +
              \left(\begin{array}{cc}f & 0 \\ 0 & g \end{array}\right)\\
  = & \left(\begin{array}{cc}f & 0 \\ y & g \end{array}\right)  +
              \left(\begin{array}{cc}f & x \\ 0 & 0 \end{array}\right) \\
  = & \left(\begin{array}{cc}f & x \\ 0 & g \end{array}\right)  +
              \left(\begin{array}{cc}f & 0 \\ y & 0 \end{array}\right)
\end{align*}
that is $| \chi( \alpha ) |\geq 5$ in $L$, which is a contradiction. So
$| M | \neq 0,  | N | \neq 0$ is not possible.

Without loss of generality we may assume that $| N | =0$. So
$R = \left(\begin{array}{cc}A & M \\ 0 & B \end{array}\right)$, clearly
$2 \leq | M | \leq 3 = \win(R)$. By \textbf{Lemma \ref{winlm0}}, $\win(A) \leq 3$. To prove that $| M | = 3$, on contrary let $M = \{ 0, x \}$. Assume $\win(A) = 2$, then there exists least one $a \in A$ such that $| \chi(a) |=2$.

\textbf{Case I:} Let $a = u_1 + e_1 = u_2 - e_2$ be two distinct weak clean expressions of $a$ in $A$, where $u_1, u_2 \in \U(A)$ and $e_1, e_2 \in \idem(A)$. Then
$e_1x = u_2x -u_1x - e_2x = -e_2x + x -x = -e_2x = e_2x$. If $e_1x=0$, then for
$\alpha = \left( \begin{array}{cc}   a & 0 \\    0 & 0 \\  \end{array} \right) $ we have
$\chi(\alpha) \supseteq \left\{ \left( \begin{array}{cc}   e_i & w \\    0 & 1 \\  \end{array} \right): i = 1, 2, w \in M \right\}$. Showing that $\win(R) \geq 4$, which is not possible.
If $e_1x = x$, then for $\alpha = \left( \begin{array}{cc}   a & 0 \\    0 & 1 \\  \end{array} \right)$, we have $\chi(\alpha) \supseteq \left\{ \left( \begin{array}{cc}   e_i & w \\    0 & 0 \\  \end{array} \right): i = 1, 2, w \in M \right\}$, showing that $\win(R) \geq 4$, which is a contradiction.\\
Similarly in \textbf{Case II}, letting $a = u_1 + e_1 = u_2 + e_2$ be two distinct weak clean and in \textbf{Case III}, let $a = u_1 - e_1 = u_2 - e_2$ be two distinct weak clean expressions of $a$ in $A$, where $u_1, u_2 \in \U(A)$ and $e_1, e_2 \in \idem(A)$, we get contradictions. Therefore $\win(A) = 1$. Similarly $\win(B) = 1$. Now by \textbf{Theorem \ref{winth2}}, we have $\win(R)=2$, a contradiction, hence $| M | =3$.

Now it remains to show that, $\win(A) = \win(B) = 1$. For $e^2=e \in A$, we have
$M= eM \oplus (1-e)M$, without loss of generality let $| eM | \neq 0$. On contrary let us assume $\win(A) > 1$. So we have $a \in A$ such that $| \chi(a) | \geq 2$, that is, we have at least two distinct weak clean expressions of $a$ in $A$.\\
\textbf{Case I:} If $a = u_1 + e_1 = u_2 - e_2$, where $u_1, u_2 \in \U(A)$ and $e_1, e_2 \in \idem(A)$ such that $e_1 \neq e_2$. Let $M = e_1M$ then for $w \in M$ and
$\alpha = \left( \begin{array}{cc}   a & 0 \\    0 & 1 \\  \end{array}\right)$ we have
$\alpha= \left( \begin{array}{cc}   u_2 & 0 \\    0 & 1 \\  \end{array} \right) -
            \left( \begin{array}{cc}   e_2 & 0 \\    0 & 1 \\  \end{array} \right) =
            \left( \begin{array}{cc}   u_1 & -w \\    0 & 1 \\  \end{array} \right) +
            \left( \begin{array}{cc}   e_1 & w \\    0 & 0 \\  \end{array} \right), $
implies $\chi(\alpha) \geq 4$, a contradiction. If $e_1M = 0$,
$\alpha = \left( \begin{array}{cc}   a & 0 \\    0 & 0 \\  \end{array}\right)$ we have
$\alpha = \left( \begin{array}{cc}   u_2 & 0 \\    0 & -1 \\  \end{array} \right) -
            \left( \begin{array}{cc}   e_2 & 0 \\    0 & 1 \\  \end{array} \right)
            \left( \begin{array}{cc}   u_1 & -w \\    0 & 1 \\  \end{array} \right) +
            \left( \begin{array}{cc}   e_1 & w \\    0 & 1 \\  \end{array} \right), $
implies $\chi(\alpha) \geq 4$, a contradiction.

Similarly in \textbf{Case II}, letting $a = u_1 + e_1 = u_2 + e_2$ be two distinct weak clean and in \textbf{Case III},  letting $a = u_1 - e_1 = u_2 - e_2$ be two distinct weak clean expressions of $a$ in $A$, where $u_1, u_2 \in \U(A)$ and $e_1, e_2 \in \idem(A)$ we get contradictions. Therefore we have $\win(A) = 1$. Similarly $\win(B) = 1$.

\pagebreak

\end{document}